\documentclass[11pt,letterpaper,reqno]{amsart}
\allowdisplaybreaks 
\belowdisplayskip
\abovedisplayskip

\numberwithin{equation}{section} 

\usepackage{pgfplots}
\pgfplotsset{compat=1.18}  
\usetikzlibrary{intersections} 
\usepgfplotslibrary{fillbetween}

\usepackage{lmodern}
\usepackage[T1]{fontenc}

\makeatletter
\@namedef{subjclassname@2020}{%
  \textup{2020} Mathematics Subject Classification}
\makeatother
\usepackage[table,xcdraw]{xcolor}
\usepackage[a4paper,top=2.5cm,bottom=2.5 cm,left=2 cm,right=2 cm]{geometry}

\usepackage[utf8]{inputenc} 
\usepackage[english]{babel} 
\usepackage{changes}
\usepackage{fancyhdr}
\usepackage{amssymb}
\usepackage[toc,page]{appendix}
\usepackage{amsmath,amssymb,mathrsfs,mathrsfs,pifont,cancel,paralist,mathbbol,bbm}
\usepackage{amsthm}
\usepackage{mathtools}
\usepackage{caption, subcaption}
\usepackage{tabularx}
\usepackage{graphicx}
\usepackage{mathrsfs}
\usepackage{pdfcomment}
\usepackage{url}
\hypersetup{hidelinks}
\usepackage[a-1b]{pdfx}
\usepackage[T1]{fontenc}
\usepackage{bm}
\usepackage{enumitem}
\usepackage{comment}
\usepackage[all]{xy}
\usepackage{tikz}
\usepackage{tikz-cd}
\usepackage{mathbbol}
\usepackage{epigraph}
\usepackage{setspace}
\usepackage{ stmaryrd }
\usepackage{nicefrac,xfrac}
\usepackage{amsmath}
\usepackage{makecell} 
\usepackage{hyperref}

\definecolor{forestgreen(web)}{rgb}{0.13, 0.55, 0.13}

\newcommand{\commentoo}[1]{}
\newcommand{\N}{\mathbb{N}}
\newcommand{\Z}{\mathbb{Z}}
\newcommand{\K}{\mathbb{K}}
\newcommand{\F}{\mathbb{F}}

\newcommand{\Q}{\mathbb{Q}}

\newcommand{\R}{\mathbb{R}}
\newcommand{\C}{\mathbb{C}}

\newcommand{\D}{\mathbb{D}}

\newcommand{\mcal}{\mathcal}

\newcommand{\ord}[1]{{\left| {#1}\right|}}

\newcommand{\jac}[2]{\left(\frac{#1}{#2}\right)}
\newcommand{\sfree}{\text{ square-free}}
\newcommand{\leg}[2]{\left(\frac{#1}{#2}\right)}
\newcommand{\summ}[1]{\sum_{\substack{#1}}}
\newcommand{\mbf}[1]{\mathbf{#1}}
\theoremstyle{definition}
\newtheorem{defi}{Definition}[section]

\newtheorem{example}[defi]{Example}

\newtheorem{remark}[defi]{Remark}

\newtheorem*{outline}{Outline of the paper}
\newtheorem*{ai}{AI disclosure statement}

\theoremstyle{plain}
\newtheorem{prop}[defi]{Proposition}

\newtheorem{lemma}[defi]{Lemma}
\newtheorem{theorem}[defi]{Theorem}

\newtheorem*{ack}{Acknowledgements}

\title{2-Selmer in arithmetic progressions}
\author{Giacomo Bortolussi}
\address{Università di Roma Tor Vergata\\ Dipartimento di Matematica\\ 00133\\ Rome}
\email{bortolussi@axp.mat.uniroma2.it}
\begin{document}
\begin{abstract}
We prove uniform asymptotics for the 2-Selmer cardinality in the congruent number family $y^2=x^3-D^2x$ with twisting parameter $D$ restricted to long arithmetic progressions.
\end{abstract}
\maketitle

\setcounter{tocdepth}{1}

\section{Introduction}

For a positive integer $D$, let $E_D$ be the congruent number curve $$y^2=x^3-D^2x.$$ Since $E_D$ is a quadratic twist of $E_1$, it is natural to restrict to square-free $D$. The rank $r(E_D)$ can be studied through the $2$-Selmer group: as $E_D(\Q)[2]\simeq(\Z/2\Z)^2$, the quantity
\[
s(E_D)\coloneqq\dim_{\F_2}\mathrm{Sel}_2(E_D/\Q)-2
\]
satisfies $r(E_D)\le s(E_D)$. In his seminal work~\cite{HB}, Heath-Brown proved that the average of $2^{s(E_D)}$ over square-free $D$ in each odd residue class modulo $8$ is $3$. He subsequently determined the higher moments of the $2$-Selmer cardinality in this family~\cite{HB2}.

We establish a uniform version of the first-moment asymptotic over long arithmetic progressions.

    \begin{theorem}
    \label{teo_HB_arithm_prog} Fix $N,\epsilon >0$. Then, for every $T>3$ and integers $q,a\le (\log T)^N$, we have $$\sum_{\substack{D\le T\sfree\\D\equiv a\bmod q}}2^{s(E_D)}=3\#\left\{D\le T\,\,\Bigg|\,\,\begin{aligned}&D\text{ square-free}\\[-1.5mm]&D\equiv a\bmod q
    \end{aligned} \right\}+O\left(\frac{T}{q(\log T)^{1/4-\epsilon}}\right)\,,$$ where the constant in the error term depends only on $N$ and $\epsilon$.
\end{theorem}

In general, both the algebraic rank and the $2$-Selmer rank exhibit parity biases, which, in the case of the congruent number curve, are determined by the class of $D\bmod 8$.
Theorem \ref{teo_HB_arithm_prog} implies that the average size of the $2$-Selmer group has no asymptotic bias in any fixed arithmetic progression, even when the modulus tends to infinity as slowly as $q\sim (\log T)^N$. Beyond this logarithmic regime, however, a bias may appear; see Example \ref{ex_sap}.

The proof contains two main novelties. First, we adapt Heath-Brown's $2$-descent counting argument to arithmetic progressions, uniformly for moduli of logarithmic size. When the residue class is not coprime to the square part of the modulus, certain character sums cease to exhibit the cancellation used in the original argument. We separate the prime factors responsible for this obstruction and exploit the combinatorics of the sixteen descent variables to find another source of cancellation. We also use Nair's estimate for multiplicative functions in arithmetic progressions to uniformly control the remaining terms. Second, we use Monsky's matrix description to extend the explicit $2$-descent counting formula to even twists.

Following Heath-Brown's work~\cite{HB,HB2}, Swinnerton-Dyer~\cite{SD} established limiting distributions for the $2$-Selmer ranks of suitable quadratic twist families with full rational $2$-torsion, ordering twists by their number of prime factors. Kane~\cite{Ka} obtained the corresponding distributions with respect to natural density. Klagsbrun, Mazur and Rubin~\cite{KMR} developed a Markov model for Selmer-rank distributions and studied parity phenomena. Bhargava--Shankar~\cite{BS}, using a structurally different method, established a similar asymptotic for the family of all elliptic curves. More recently, Smith~\cite{Sm} determined the distribution of higher $2$-power Selmer groups in suitable twist families, in parallel with his results on the $2$-primary parts of quadratic class groups.

There is a close analogy between the statistics of the $2$-Selmer rank in the congruent number family and those of the $4$-rank of quadratic class groups, the latter being related to the Cohen–Lenstra heuristics. In the class-group setting, Fouvry and Klüners~\cite{FK} established results analogous to those of Heath-Brown for $2$-Selmer groups. Subsequently, Milovic~\cite{Mi} extended the work of Fouvry--Klüners to arithmetic progressions with fixed square-free modulus and coprime residue class, adapting techniques from Heath-Brown's method. Our result provides the corresponding extension of Heath-Brown's theorem in the context of $2$-Selmer groups, with two key differences: uniformity for moduli of logarithmic size and the removal of the coprimality assumption on the residue class.
The latter, in particular, requires isolating and controlling certain error terms that obstruct character cancellation in arithmetic progressions.

\begin{outline}
Section \ref{sec_prel} recalls the $2$-descent background, establishes the counting formula for even twists, and describes the combinatorics of linked descent variables. In Section \ref{sec_2sap} we prove the uniform asymptotic for odd twists by separating the various types of error terms, and then treat even twists by adapting the main-term analysis. We conclude with a counterexample for short arithmetic progressions.
\end{outline}

\begin{ai}
GPT-5.6 Sol and GPT-6 Astra were used under human supervision to prove Lemma \ref{lemma_struct} and to run calculations for the expression for $g_2$ in \eqref{defg_even}. Every other idea and mathematical argument in this paper is due to the author.
\end{ai}
\begin{ack}
    The author would like to thank his supervisor, Efthymios Sofos, for introducing him to the topic and for the support during the writing process. 
\end{ack}

\section{Preliminaries}\label{sec_prel}

\subsection{Algebraic background}
Given an elliptic curve $E$ defined over $\Q$, Mordell's theorem asserts that $$E(\Q)\cong E_{\mathrm{tors}}\oplus \Z^{r(E)}\,.$$
While the torsion part is well understood, the algebraic rank (over $\Q$) $r(E)$ is generally hard to access. The simplest way to control the algebraic rank is through Selmer groups: one can look at the embedding $$\frac{E(\Q)}{2E(\Q)}\hookrightarrow\mathrm{Sel}_2(E/\Q)\,,$$ and then $2^{r(E)}\cdot \#E(\Q)[2]\le \#\mathrm{Sel}_2(E/\Q)$. One defines $s(E)\coloneqq \dim_{\F_2}\left(\mathrm{Sel}_2(E/\Q)\right)-\dim_{\F_2}\left(E(\Q)[2]\right)$; then $r(E)\le s(E)$. This whole theory fits into a more natural cohomological framework~\cite{Si}, but here, for simplicity, we give only the strictly necessary notions. Each element of the $2$-Selmer group of $E$ can be viewed as an everywhere locally solvable system of two quadrics in $\mathbb P^3$ (see~\cite{HB}), called a 2-descent, and each globally solvable system corresponds to a point in $E(\Q)/2E(\Q)$. In this sense, the cokernel of the map above (which is the 2-torsion part of the Tate-Shafarevich group) represents the obstruction to the local-global principle.
For the congruent number family $E_D$, after quotienting by the image of rational $2$-torsion, for $D $ odd and square-free the descent classes are represented by systems \begin{equation}\label{system}\begin{cases}
    D_1X^2+D_4W^2=D_2Y^2\\D_1X^2-D_4W^2=D_3Z^2
\end{cases},\end{equation} where $D=D_1D_2D_3D_4$ and $D_i\in\Z_{>0}$ for $1\le i\le4$. There are $2^{r(E_D)}$ globally solvable systems and $2^{s(E_D)}$ everywhere locally solvable systems of this form~\cite[Lemma 1]{HB}.
The system (\ref{system}) is always solvable over $\R$, and the local solvability needs to be checked only at the places dividing $2D$. If we write $s(D)\coloneqq s(E_D)$ and $r(D)\coloneqq r(E_D)$, then $s(D)\leq 2\omega(D)$ and
\begin{equation}
\label{bound-sel}
2^{s(D)}\cdot 4^{\omega(D)}\,.
\end{equation}

\subsection{Counting 2-descents} Using Jacobi symbols, Heath-Brown~\cite[Lemma 3]{HB} produced an indicator function that detects whether any system of the form (\ref{system})  is everywhere locally solvable, and this led to the $2$-descent counting formula (\ref{2desc}). This approach was later generalized by Monsky~\cite[Appendix]{HB2}.

Let $$\mcal I\coloneqq \Big\{(u,v): u\in \{1,\ldots,4\},\, v\in \{0,\ldots,4\},\, v\neq u\Big\} $$ be a set of indices and let $\{D_{uv}\}_{(u,v)\in \mcal I}$ be 16 integer variables indexed by $\mcal I$. For a fixed $k\in \Z_{>0}$ let $$\D_k(D)=\{(d_1,\ldots,d_k)\in\mathbb Z_{>0}^k:d_1\cdot\ldots \cdot d_k=D\}\subset\N^{k}\,.$$ We identify the set of all possible decompositions $D=\prod_{(u,v)\in \mcal I} D_{uv}$ with $\D_{16}(D)$. Since we will be working with odd square-free integers $D$, in any decomposition $ \mathbb D_{16}(D)\ni\mcal D=\{D_{uv}\}$ all $D_{uv}$ will be square-free and pairwise coprime. For two different indices $(i,j), (k,l)\in\mcal I$, we say that $(i,j)\preceq(k,l)$ if $j\neq0$ and $k\notin\{i,j\}$ (note that $\preceq$ is not an order). For any odd, square-free $D$ and for any $\mcal D\in \D_{16}(D)$, define \begin{equation}\label{defg} g(\mcal D)\coloneqq4^{-\omega(D)} \leg{-1}{\alpha(\mcal D)}\leg{2}{\beta(\mcal D)}\prod_{(i,j)\preceq (k,l)}\left(\frac{D_{kl}}{D_{ij}}\right)\,,
\end{equation}  
where $\alpha(\mcal D)\coloneqq D_{12}D_{14}D_{23}D_{21}$, $\beta(\mcal D)\coloneqq D_{24}D_{31}D_{34}D_{21}$. Then, by explicitly counting 2-descents on $E_D$, we have \begin{equation}\label{2desc}
    \frac 1 4\#\mathrm{Sel}_2(E_D/\Q)=2^{s(D)}=\sum_{\mcal D\in \mathbb D_{16}(D)} g(\mcal  D)\,.\end{equation}

In particular, by bounding all the Jacobi symbols in $g$ by $1$, one recovers $2^{s(D)}\le 4^{\omega(D)}$.

\subsection{Linked and joined variables} 
Formula (\ref{defg}) is influenced by quadratic reciprocity; therefore, it is important to keep track, for every pair of variables $D_{ij},D_{kl}$, of whether both, exactly one, or neither of $\leg {D_{ij}}{D_{kl}}$ and $\leg{D_{kl}}{D_{ij}}$ appears in (\ref{defg}).
Following the terminology introduced by Heath-Brown, we say that two indices $(i,j),(k,l)\in\mcal I$ are:\begin{itemize}
    \item \textbf{joined} if both $(i,j)\preceq (k,l)$ and $(k,l)\preceq(i,j)$. This means that the Jacobi symbols $\jac{D_{ij}}{D_{kl}}$ and $\jac{D_{kl}}{D_{ij}}$ appear in (\ref{defg}), and their contribution to (\ref{defg}) reduces, by quadratic reciprocity, to $(-1)^{\frac{D_{ij}-1}{2}\frac{D_{kl}-1}{2}}$. Equivalently, we ask $(i\neq k)\wedge(j\neq 0,i,k)\wedge (l\neq 0,i,k)$;
    \item \textbf{independent} if they are not joined;
    \item \textbf{linked} if exactly one of $(i,j)\preceq (k,l)$ and $(k,l)\preceq(i,j)$ holds. That is, exactly one of the Jacobi symbols $\jac{D_{ij}}{D_{kl}}$ and $\jac{D_{kl}}{D_{ij}}$ appears in (\ref{defg}). Equivalently, we ask $(i\neq k)\wedge\left[(j\neq 0,k)\oplus (l\neq 0,i)\right]$;
    \item \textbf{unlinked} if they are not linked.
\end{itemize}

The following lemma clarifies the structure underlying linked variables.
\begin{lemma}\label{lemma_struct}
    Rearrange the 16 indices in $\mcal I$ as 
    \begin{equation}\label{table_indexes}
\begin{array}{cccc}
10&34&42&23\\
43&20&14&31\\
24&41&30&12\\
32&13&21&40
\end{array}.
\end{equation}
    Then two indices are linked if and only if they belong to the same row or column.
\end{lemma}
\begin{proof} In~\cite{HB2} Heath-Brown identifies  $\mcal I$  with $\F_2^4$ and relabelled  $\{D_{ij}\}_{(i,j)\in \mcal I}$ as $\{D_{x,y,u,v}\}_{(x,y,u,v)\in \F^4_2}$ according to Table (\ref{table_HB}). In this way, two indices $(x_1,y_1,u_1,v_1),(x_2,y_2,u_2,v_2)$ are linked if and only if    $$Q\big((x_1,y_1,u_1,v_1)+(x_2,y_2,u_2,v_2)\big)\coloneqq(x_1+x_2)(v_1+v_2)+(y_1+y_2)(u_1+u_2)$$ is $1$.
Changing the basis of the index space according to $$(a,b,c,d)=\Phi(x,y,u,v)
=(x+u+v, y+u, x+v,x+y+u).
 $$ one can relabel the variables as $\overline D_{a,b,c,d}=D_{\Phi^{-1}(a,b,c,d)}$, and in the new basis $Q$ becomes   \begin{align*}\overline Q\big( (a_1,b_1,c_1,d_1),(a_2,b_2,c_2,d_2)\big)&= \big((1+a_1+a_2)(1+b_1+b_2)\big)+\big((1+c_1+c_2)(1+d_1+d_2)\big)\\
&= \big((a_1=a_2)\wedge(b_1=b_2)\big)\oplus \big((c_1=c_2)\wedge(d_1=d_2)\big)\,,\end{align*}
where, from the first to the second line, we translated the expression from $(\F_2,+,\cdot)$ into the same expression under logical operators.
From the expression above, one can create a $4\times 4$ table with the variables $\overline D_{abcd}$, where the rows are indexed by $\{(a,b)\}_{a,b\in \F_2}$ and the columns are indexed by $\{(c,d)\}_{c,d\in \F_2}$; in this way, two variables are linked if and only if they are in the same row or in the same column. Going back to our original labelling of the variables, one obtains, for example, the table in the statement.
\end{proof}
This reformulates and strengthens the combinatorial description  of~\cite[Lemma 9]{HB}, as now 4-tuples of pairwise unlinked variables are generalized diagonals in the table (\ref{table_indexes}).

\subsection{2-descent counting for even twists} Due to the irreguar behavior of squares in $\Q_2$, formula (\ref{2desc}) does not apply to even twists.
Heath-Brown~\cite[page 5]{HB2} stated that \lq\lq \textit{One can do an analogous calculation for even values of $D$, again based on the
method of Monsky’s appendix}\rq\rq. Already Chan~\cite{Ch} analyzed $2^{s(2D)}$, but for her purposes she dropped the local solubility condition in $\Q_2$, obtaining only an upper bound.
Here, we extend the $2$-descent counting formula to even values of the twisting parameter, building on Monsky's matrix approach.

For any odd, square-free $D$ and for any $\mcal D\in \D_{16}(D)$, define \begin{equation}\label{defg_even} g_2(\mcal D)\coloneqq\leg{-1}{\alpha_2}\leg{2}{\beta_2}4^{-\omega(D)}\prod_{(i,j)\preceq (k,l)}\left(\frac{D_{kl{}}}{D_{ij}}\right)\,,
\end{equation}  
where $\alpha_2=D_{30}D_{31}D_{32}D_{34}$ and $
\beta_2=D_{13}D_{23}D_{31}D_{32}D_{34}D_{43}$. We prove that, by explicitly counting 2-descents on $E_{2D}$, we have \begin{equation}\label{2desc_even}
    \frac 1 4\#\mathrm{Sel}^2(E_{2D})=2^{s(2D)}=\sum_{\mcal D\in \mathbb D_{16}(D)} g_2(\mcal  D)\,.\end{equation}
As a direct consequence, we obtain 
\begin{equation}
    \label{bound_sel_even}{2^{s(2D)}}\le 4^{\omega(D)}\,.
\end{equation}
For $p$ an odd prime and $a\in \Z$ coprime to $p$, define $[a/p]$ to be the additive Legendre symbol, i.e. the unique element of $\F_2$ such that $(-1)^{[a/p]}=\leg a p$. If $D=p_1\cdots p_k$ with $p_r$ distinct primes and $k=\omega(D)$, define
$$a_{rs}\coloneqq [p_s/p_r]\text{ for }r\neq s,\qquad a_{rr}\coloneqq \sum_{\substack{s=1,\ldots,k\\s\neq r}}a_{r,s}\,,\qquad \mathbf A\coloneqq\{a_{rs}\}_{1\le r,s\le k}$$ 
and, for $u\in \{-1,2\}$, 
$$\mathbf D_u\coloneqq \mathrm{diag}([u/p_1],\ldots,[u/p_k])\,.$$

By the work of Monsky~\cite[Appendix]{HB2}, we know that 
\begin{equation*}2^{s(2D)}=\#\ker M_{2D},\,\qquad \text{where}\qquad M_{2D}=
\begin{pmatrix}
\mathbf D_2&\mathbf A+\mathbf D_2\\
\mathbf A^T+\mathbf D_2&\mathbf D_{-1}
\end{pmatrix}\,.
\end{equation*}

Given $(x,y)\in \mathbf F_2^{2k}$, by character orthogonality, we have
$$\sum_{u,v\in\mathbf F_2^k}(-1)^{(u,v)M_{2D}(x,y)^T}
=
\begin{cases}4^{k}&\text{ if }M_{2D}(x,y)^T=0\,,\\
0&\text{ if } M_{2D}(x,y)^T\ne0\,.\end{cases}$$
Consequently,
\begin{equation}\label{2desc.2}2^{s(2D)}=4^{-\omega(D)}\sum_{x,y,u,v\in\mathbf F_2^{k}}(-1)^{Q(x,y,u,v)}\coloneqq \sum_{x,y,u,v\in \F_2^{k}}g_2(x,y,u,v)\end{equation}
where
\begin{align*}Q(x,y,u,v)
=&
u\mathbf D_2x^T+u(\mathbf A+\mathbf D_2)y^T+v(\mathbf A^T+\mathbf D_2)x^T+v\mathbf D_{-1}y^T\\
=& \sum_{r=1}^k [-1/p_r] v_ry_r 
\,+ \,\sum_{r=1}^k [2/p_r] \bigl(u_rx_r+u_ry_r+v_rx_r\bigr)
\\+&\,
\sum_{\substack{1\le r,s\le k\\r\ne s}} [p_s/p_r]\bigl(u_r(y_r+y_s)+x_r(v_r+v_s)
\bigr),
\end{align*}
Observe that the non-diagonal part of the formula above is the same for both even and odd twists. Since $g_2(\cdot)=(-1)^{Q(\cdot)}$, we have
\begin{equation}\label{defg2.2}
    g_2(x,y,u,v)=4^{-\omega(D)}\prod_{\substack{r=1,\ldots, k\\v_ry_r=1}}\leg {-1} {p_r}\, \prod_{\substack{r=1,\ldots, k\\u_rx_r+u_ry_r+v_rx_r=1}}\leg {2} {p_r}\,  \prod_{\substack{r,s=1,\ldots, k, \,r\neq s\\u_r(y_r+y_s)+x_r(v_r+v_s)=1}}\leg {p_s} {p_r}\,.\end{equation}

Given $(x,y,u,v)\in (\F_2^k)^4$, let $\sigma_r(x,y,u,v)=(x_r,y_r,u_r,v_r)\in \F^4_2$. For any $\sigma\in \F^4_2$, define $$D_\sigma(x,y,u,v)=\prod_{\sigma_r(x,y,u,v)=\sigma}p_r\,.$$ Then $\prod_{\sigma\in \F^4_2}D_\sigma(x,y,u,v)=D$ for every $(x,y,u,v)\in (\F_2^k)^4$.
Relabel the $D_\sigma$ according to Heath--Brown's sixteen indices $\{(i,j)\}_{\mcal I}$, using
\begin{equation}\label{table_HB}
\begin{array}{c|cccc}
i=1&
D_{10}=D_{0001}&
D_{12}=D_{1011}&
D_{13}=D_{1001}&
D_{14}=D_{0011}\\[1mm]
i=2&
D_{20}=D_{0100}&
D_{21}=D_{1110}&
D_{23}=D_{0110}&
D_{24}=D_{1100}\\[1mm]
i=3&
D_{30}=D_{0101}&
D_{31}=D_{1101}&
D_{32}=D_{0111}&
D_{34}=D_{1111}\\[1mm]
i=4&
D_{40}=D_{0000}&
D_{41}=D_{0010}&
D_{42}=D_{1000}&
D_{43}=D_{1010}.
\end{array}
\end{equation}

With this change of variables, one checks that the first product in (\ref{defg2.2}) gives $\leg {-1}{\alpha_2}$, the second product gives $\leg 2{\beta_2}$, and the third product becomes the product of $\leg{D_{kl}}{D_{ij}}$ with the condition $k\neq i,j$. Finally, (\ref{defg2.2}) and (\ref{2desc.2}) become, respectively, (\ref{defg_even}) and (\ref{2desc_even}).

\section{2-Selmer on long arithmetic progressions} \label{sec_2sap}

Let $q,a\in\Z_{>0}$ with $q,a\le (\log T)^N$, and consider the arithmetic progression $D\equiv a\bmod q$. Let $$S(T;q,a)\coloneqq \summ{D\le T\\D\equiv a\bmod q\\D\sfree }2^{s(D)}\,.$$ Let $h^2$ be the  largest square divisor of $q$, i.e. the unique positive integer $\le (\log T)^{N}$ such that $h^2\mid q$ and $q/h^2$ is square-free.

Milovic~\cite{Mi} generalized the work of Fouvry--Klüners~\cite{FK} on the statistics of $\#2\mathrm{Cl}(\Q(\sqrt D ))[4]$ by restricting $D$ to arithmetic progressions $\{qn+a\}$ with $\gcd(a,q)=1$. His method extends to progressions with $\gcd(a,h)=1$, but when $\gcd(a,h)\neq 1$ an obstruction to character cancellation appears in what he calls (C4) multi-intervals. The same type of obstruction arises for 2-Selmer groups when one attempts to impose congruence conditions in Heath-Brown's argument (more precisely, the obstruction lies in the application of the Siegel--Walfisz method~\cite[Lemma 6]{HB}). 

To overcome this issue, we modify Heath-Brown's error-term analysis to isolate the obstruction to character cancellation, and we handle the problematic terms in a different way.
In doing so, under the extra assumption $q\le (\log T)^N$, we also achieve uniformity of the error term in $q,a$.

\begin{theorem}
    \label{teo_HB_arithm_prog_odd_restr} Fix $N,\epsilon>0$. Let $T>3$ and $h\in\Z_{>0}$ with $h\le(\log T)^{N/2}$. Then, for every $q,a\in\Z_{>0}$ with $q,a\le (\log T)^N$, $q$ even, $a$ odd, $h^2\mid q$, and $\mu^2(q/h^2)=1$, we have $$S(T;q,a)=3\#\left\{D\le T\bigg|\,\begin{aligned}
        &D\equiv a\bmod q\\[-1.5mm]&D\sfree
    \end{aligned}\right\}+O\left(\frac T{q\,(\log T)^{1/4-\epsilon}}\right)\,,$$
    where the implied constant depends only on $N,\epsilon$.
\end{theorem}
\begin{remark}
        In the literature, there are two main approaches to the study of $2$-Selmer statistics in quadratic twist families. The first, pioneered by Heath-Brown~\cite{HB,HB2} and adopted in the present work, orders the twisting parameter $D$ by height, typically averaging over $D\leq X$. The second, due to Swinnerton-Dyer~\cite{SD}, regroups the twists by the value of $\omega(D)$ and studies the limit distribution under $\omega(D)\to \infty$. For fixed $\omega(D)$ we also mention the work of Rhoades~\cite{Ro}.  Kane~\cite{Ka} established the compatibility of these two regimes, showing that the Swinnerton-Dyer distribution also arises under the natural ordering by height. Having $\gcd(q,a)\neq 1$ establishes a sort of mixed regime, where $D$ varies linearly but also has a fixed set of divisors.
\end{remark}
\begin{remark}
    Even though the error term in Theorem \ref{teo_HB_arithm_prog_odd_restr} does not depend on $h$, the combinatorial separation of the implied error terms does, so we preferred to keep that dependence explicit.
\end{remark}

First, observe that if $\gcd(q,a)$ is not square-free, then the whole sum is empty. Thus, we may safely assume that $\gcd(q,a)$ is square-free.
Let $\ell \coloneqq \gcd(D,h)$, and decompose any square-free $D\in\Z_{>0}$ as $D=\ell K$. Note that $\ell \mid \gcd(a,q)$, thus $\ell$ must also be square-free.
\begin{remark}The choice of $\ell$ depends on $D$ and $h$, and, even though $D\equiv a\bmod q$ forces $\ell=\gcd(a,h)$, the choice of $\ell$ is in principle made depending only on $D$ and $h$. \end{remark}
 Call $q'=q/\ell,\, a'=a/\ell$. For  $D\le T$ respecting $D\equiv a\bmod q$, we have $K\le T/\ell$, $K\equiv a'\bmod q'$, $\mu^2(K)=1$, and $\gcd(K,\ell)=1$.
For any $\mcal D\in \D_{16}(D)$ write \begin{align*}\ell_{uv}=\gcd(D_{uv},\ell),\qquad K_{uv}=D_{uv}/\ell_{uv}\qquad \text{and  }\qquad\qquad \mcal D=\mcal L\cdot \mcal K\,,\qquad\text{where}\\
\qquad \D_{16}(\ell)\ni \mcal L=\{\ell_{uv}\}_{(u,v)\in \mcal I}\qquad\text{ and }\qquad\D_{16}(K)\ni \mcal K=\{K_{uv}\}_{(u,v)\in \mcal I}\,.\qquad\end{align*}
Observe that, since $D$ is odd and square-free, all $\ell_{uv}\,,K_{uv}$ must be odd, square-free, and pairwise coprime. 

If one tries to apply verbatim Heath-Brown's method to arithmetic progressions, when $\ell\neq 1$ characters' orthogonality may transform some non-trivial character into a  trivial character, interfering then with  the application of Siegel--Walfisz. Splitting $\mcal D$ as $\mcal L\mcal K$ allows to separate the obstruction to characters' cancellation $\mcal L$ to the source of characters' cancellation $\mcal K$.
 We then split the contribution of $\mcal L\mcal K$ to $g$ as $$g(\mcal L\mcal K)=g(\mcal L)g(\mcal K)R(\mcal L,\mcal K)\qquad\text{ with }\qquad R(\mcal L,\mcal K)=\prod_{(i,j)\preceq(k,l)}\leg{\ell_{k,l}}{K_{ij}}\leg{K_{kl}}{{\ell_{ij}}}\,,$$ then $$2^{s(D)}=2^{s(\ell K)}=\sum_{\mcal L\in \D_{16}(\ell)}g(\mcal L)\sum_{\mcal K\in \D_{16}(K)}g(\mcal K)R(\mcal L,\mcal K)\,.$$

In the first part of our argument, we show that Heath-Brown's method can be carried out by isolating the problematic values of $\mcal L$. In the second part, we use the combinatorics of linked indices to show that the contribution of certain problematic choices of $\mcal L$ is negligible. In the final part, we reunite $\mcal L\mcal K$ into $\mcal D$ and conclude by Heath-Brown's argument. 

For $B,\eta>0$ (to be chosen later), set \begin{equation}\label{def_zC}z(T)\coloneqq2^{\lfloor{\log_2(\log T)^B}\rfloor}\qquad \text{and}\qquad\qquad C(T)\coloneqq 2^{\lfloor(\log T)^\eta/\log 2\rfloor}\,.\end{equation} Then $
\tfrac12(\log T)^B<z\le(\log T)^B$ and $\tfrac12e^{(\log T)^\eta}<C\le e^{(\log T)^\eta}.
$ Note that Heath-Brown's original choice of $C$ is incorrect, we adopt Fouvry--Kluner's correction \cite[Definition (36)]{FK}.
We now split the set $$\left\{\mcal D\in \D_{16}(D)\bigg|\begin{aligned}
    &D\le T\sfree \\[-1.5mm]
    &D\equiv a\bmod q
\end{aligned}\right\}$$
into different subsets, each contributing either to a main term or to an error term. In particular: \begin{itemize}
    \item Type 1 and R error terms rely on Nair's bound\cite{Na};
    \item Type 2 error terms rely on Wilson's hyperbola method for Jacobi symbols\cite{Wi};
    \item Type 3 and 4 error terms rely, in different ways, on Siegel--Walfisz' method.
\end{itemize} We remark that our subdivision depends on the choices of $h,B,\eta$.

\subsection{Type 1 error terms} We say that $\mcal D\in \D_{16}(D)$ is a type 1 term  if $\mcal D$ contains at most three variables $D_{ij}$ larger than $C(T)$. Let
\begin{align*}&G^{(1)}(D)\coloneqq\mu^2(D)\summ{\mcal D\in \D_{16}(D)\\\mcal D\text{ of type 1}}g(\mcal D)\qquad\text{and}\qquad  S^{(1)}(T;q,a)\coloneqq \sum_{\substack{D\le T\\D\equiv a\bmod q\\D\sfree}}G^{(1)} ( D)\,.
\end{align*}
To deal with the type 1 terms, we employ Nair's lemma~\cite{Na}.

\begin{lemma}[Nair]\label{lemma_nair}
    Let $x,y>1$, $q,a\in \N^+$ with $\gcd(q,a)=1$ and $u>0$. Assume that for some $\alpha,\delta\in (0,1)$ and $c=c(\delta,\alpha)>0$ we have $x^{\alpha}\le y \le x$ and $q,a\le    cx^{\delta}$. Then $$ \sum_{x-y\le n\le x}u^{\omega(qn+a)}\ll C_{u,q}y (\log x)^{u-1}\,, $$ where the implied constant depends only on $u,\delta,\alpha$ and $C_{u,q}=\prod_{p\mid q}\frac{e^{-u/p}} {1-1/p}$.
\end{lemma}
\begin{remark}
    Uniform asymptotic results that are available for multiplicative functions over long arithmetic progressions often cannot be proved over short arithmetic progressions. Nair's work, on the other hand, guarantees a uniform upper bound over short arithmetic progressions for suitable multiplicative functions.
\end{remark}

\begin{prop}\label{prop_error1}
    In the setting of Theorem \ref{teo_HB_arithm_prog_odd_restr}, $0<\eta<1/16  $ and let  $C(T)$ be as in (\ref{def_zC}). Then $$\ord{S^{(1)}(T;q,a)}\ll \frac{T}{q(\log T)^{1/4-4\eta}}\,,$$ where the implied constant depends on $N,\eta$.
\end{prop}
\begin{proof}
    By bounding $\ord{g(\mcal D)}\le 4^{-\omega(D)}$, we have
\begin{align*}
    \ord{S^{(1)}(T;q,a)}&\le\summ{D\le T\\D\equiv a\bmod q}4^{-\omega(D)}\mu^2(D)\#\{\mcal D\in \D_{16}(D)| \text{ at most three }D_{ij}\ge C\}\end{align*}
    Split $\mcal D=\mcal D'\times\{D_1\}\times\{D_2\}\times\{D_3\}$ into $\binom{16}{3}$ ways, force every variable $D_{ij}$ in $\mcal D'$ to be smaller than $C$, and distribute the remaining factors of $n\coloneqq D/D'$ among $D_1,D_2,D_3$ in $3^{\omega(n)}$ ways. This counting overcounts the set above, but since all terms are nonnegative, this is acceptable. Then, \begin{equation}\label{bound_g3}G^{(1)}(D)\ll\mu^2(D)\summ{D'\le C^{13}\\D'\mid D}\leg {13}4^{\omega(D')}\leg 3 4 ^{\omega(D/D')}\end{equation} and 
    \begin{align*}
    \ord{S^{(1)}(T;q,a)}\ll\summ{D'\le C^{13}}\frac{\mu^2(D')\tau_{13}(D')}{4^{\omega(D')}}\summ{ n\le T/D'\\\gcd(n,D')=1\\nD'\equiv a\bmod q}\mu^2(n)\leg 3 4 ^{\omega(n)}\,.
\end{align*}
One cannot use Lemma \ref{lemma_nair} directly, as $a$ and $q$ might not be coprime. Letting $u=\gcd(a,q)$, $a''\coloneqq a/u$, $q''\coloneqq q/u$, write $u=\lambda\kappa$ so that $\lambda=\gcd(u, n)$, $\kappa=\gcd(u,D')$, and let $n_1=n/\lambda$, $D_1'=D'/\kappa$.  Then $\gcd(a'',q'')=1$,  and $D\equiv a\bmod q$, forces $D_1'n_1\equiv a''\bmod q''$, and so $D_1',n_1$ must be coprime to $q''$. The quantity above rewrites as
\begin{equation*}
    \summ{u=\lambda\kappa}\frac{3^{\omega(\lambda)}13^{\omega(\kappa)}}{4^{\omega(\lambda\kappa)}}\summ{D_1'\le C^{13}/ \kappa\\\gcd(D'_1,\lambda\kappa q'')=1}\frac{\mu^2(D_1')13^{\omega(D_1')}}{4^{\omega(D'_1)}}\summ{ n_1\le T/ D_1'\kappa\lambda\\\gcd(n_1, \lambda\kappa D_1q'')=1\\n_1\equiv a''/D'_1\bmod q'' }\mu^2(n_1)\leg 3 4^{\omega(n_1)}\,.
\end{equation*}
Considering the sum over $n_1$ in the expression above, we may remove the coprimality and square-free conditions on $n_1$ and apply Lemma \ref{lemma_nair} (which is allowed since $\gcd(a''/D_1',q'')=1$) to obtain a bound of $\ll \mcal C_{3/4,q''}\, T/\bigl((\log T)^{1/4}q''\lambda\kappa D_1\bigr)$, where $\mcal C_{3/4,q''}=\prod_{p\mid q''}\frac{e^{-3/(4p)}}{1-1/p}\ll2^{\omega(q'')}$. Then the whole sum can be bounded by
\begin{align*}
    &\ll \frac{T2^{\omega(q'')}}{(\log T)^{1/4}q''} \summ{u=\lambda\kappa}\frac{3^{\omega(\lambda)}13^{\omega(\kappa)}}{4^{\omega(\lambda\kappa)}\lambda\kappa}\summ{D_1'\le C^{13}/ \kappa\\\mcal \gcd(D'_1,\lambda\kappa q'')=1} \frac{\mu^2(D_1')13^{\omega(D_1')}}{4^{\omega(D'_1)}D_1'}
    \ll\frac{T 2^{\omega(q'')} 4^{\omega( u)}\log\left(C^{13}\right)^{13/4}}{q''u(\log T)^{1/4}}\,.
\end{align*}
Recalling that $q=q'' u\le (\log T)^{N}$ and $C(T)\asymp e^{(\log T)^\eta}$, we have $2^{\omega(q'')}4^{\omega(u)}\ll_{N,\eta} (\log T)^{3\eta/4} $ and we can conclude.
\end{proof}

\textbf{Dyadic splitting.}
Given $\{A_{ij}\}_{(i,j)\in \mcal I}\in \{2^m\mid m\in \Z_{\ge 0}\}^{16}$, define the corresponding \textit{dyadic multi-interval} as $$\textbf{A}=\prod_{(i,j)\in \mcal I}[A_{ij},2A_{ij})\,.$$ Observe that, restricting to $D\le T$, there are $O((\log T)^{16})$ relevant dyadic multi-intervals. In what follows, we regroup the values of $\mcal D$ by dyadic multi-intervals. We say that a multi-interval $\mbf A$ is of type 1 if at most three indices $(u,v)\in \mcal I$ satisfy $A_{uv}\ge C$.

\subsection{Type 2 error terms} We say that a dyadic multi-interval $\textbf{A}$ is of \textbf{type 2} if it is not of type 1, if there exist two linked indices $(i,j),(k,l)\in \mcal I$ such that $A_{ij},A_{kl}\ge z$. We further say that $\mcal D$ is a type 2 term if $\mcal D$ belongs to a type 2 multi-interval $\textbf{A}$.

For each type 2 multi-interval $\mathbf A$ and each odd $D$, define \begin{align*}G^{\mathbf A}(D)\coloneqq \mu^2(D)\sum_{\mcal D\in \D_{16}(D)}g(\mcal D)\mathbb 1_{\mathbf A}(\mcal D)\qquad \text{and}\qquad S^{\mathbf A}(T;q,a)\coloneqq \sum_{\substack{D\le T\\D\equiv a\bmod q\\D\sfree}}G^\mathbf A(D)
\end{align*}

Having two linked variables $D_{ij},D_{kl}$, one can, interchanging the two indices if necessary,  isolate the corresponding Jacobi symbol in $g$ as
$$g(\mcal D)=4^{-\omega(D)}b(D_{ij})c(D_{kl})\leg{D_{ij}}{D_{kl}}\,,$$ where $\ord b,\ord c\le 1$, $b$ depends on every $D_{uv}$ except $D_{kl}$ and $c$ depends on every $D_{uv}$ except $D_{ij}$.
To deal with type 2 error terms, we employ a version of the hyperbola method for Jacobi symbols due to Wilson~\cite[Theorem 1.2]{Wi}.
\begin{lemma}[Wilson]\label{lemma_hyp}
    Let $T,z\ge2$ and let $\{b_n\},\{c_m\}$ be two complex sequences with $\ord{b_n},\ord{c_m}\le 1$. Then $$\sum_{\substack{m\cdot n\le T\\z\le m,n\le T}}\mu^2(2mn)b_nc_m\jac{n}{m}\ll \frac{T\log T}{z^{1/4}}\,, $$
    where the implied constant is independent of $\{b_n\},\{c_m\}$.
\end{lemma}

\begin{prop}\label{prop_error2}
    In the setting of Theorem \ref{teo_HB_arithm_prog_odd_restr}, fix $z>2$ and let $\mathbf A$ be a type 2 multi-interval. Then $$\ord{S^\mathbf A(T;q,a)}\ll \frac{q^2T(\log T)^{9/2}}{z^{1/4}}\,,$$ where the implied constant is absolute.
\end{prop}
\begin{proof}
    Call $(i,j)$ and $(k,l)$ the two linked indices for which $A_{ij},A_{kl}\ge z$, and, without loss of generality, suppose that $\leg{D_{ij}}{D_{kl}}$ appears in the expression of $g$. Set $D=D'D_{ij}D_{kl}$ and let $\mcal D'\in \D_{14}(D')$ be such that $\mcal D=\mcal D'\times\{D_1\}\times\{D_2\}$.
    Then we have  \begin{align*}
        \ord{S^{\mathbf{A}}(T;q,a)}=\Bigg|\sum_{\substack{D\le T\\D\equiv a\bmod q}}\frac{\mu^2(D)}{4^{\omega(D)}}\sum_{\substack{\mcal D\in \D_{16}(D)\\\mcal D\in \mathbf{A}}}b_{D_1}(\mcal D')c_{D_2}(\mcal D')\leg{D_1}{D_2}{}\Bigg|
    \end{align*}
    for some choice of complex sequences $\{b_{D_1}\}, \{c_{D_2}\}$ depending on $\mcal D'$ such that $b_{D_1}$ is independent of $D_2$ and $c_{D_2}$ is independent of $D_1$.
    To simplify the notation, define \begin{align*}&b'_{D_1}(\mathcal D')=4^{-\omega(D_1)}\,\mathbb 1_{\mcal D'\in \mbf A'}\,\,  \mathbb1_{D_1\in [A_{1},2A_{1})}\,\,\mathbb 1_{(D_1,D')=1} \, b_{D_1}(\mathcal D')\,,\\
    &c'_{D_2}(\mathcal D')=4^{-\omega(D_2)} \,1_{\mcal D'\in \mbf A'}\,\, \mathbb1_{D_2\in [A_{2},2A_{2})}\,\,\mathbb1_{(D_2,D')=1}\,  c_{D_2}(\mathcal D')\,.\end{align*} In this way, \begin{align*}
        \ord{S^{\mathbf{A}}(T;q,a)}
        &\le
        \Bigg|\sum_{\substack{D\le T\\D\equiv a\bmod q}}\frac{\mu^2(D)}{4^{\omega(D)}}\sum_{\substack{\mcal D\in \D_{16}(D)\\\mcal D\in \mathbf{A}}}b_{D_1}c_{D_2}\leg{D_1}{D_2}\Bigg|
        \\&\le
        \sum_{\substack{D'\le T/ z^2\\\mcal D'\in \D_{14}(D')\\}}\frac{\mu^2(D')}{4^{\omega(D')}}\Bigg|\sum_{\substack{z\le D_1,D_2\le T/D'\\ D_1D_2\le T/D'\\D_1,D_2\text{ odd,}\sfree\\\gcd(D_1,D_2)=1\\D_1D_2D'\equiv a\bmod q}}b'_{D_1}c'_{D_2}\leg{D_1}{D_2}\Bigg|
    \end{align*}
    Bound the inner sum above by splitting it according to the values of $D_1$ and $D_2$ modulo $q$. In this way, we have $$\le \summ{a_1,a_2\bmod q\\a_1a_2D'\equiv a\bmod q}\Bigg|\sum_{\substack{z\le D_1,D_2\le T/D'\\ D_1D_2\le T/D'\\D_1,D_2\text{ odd,}\sfree\\\gcd(D_1,D_2)=1}}b''_{D_1}c''_{D_2}\leg{D_1}{D_2}{}\Bigg|\,,$$ where $b''_{D_1}=b'_{D_1}\mathbb 1_{D_1\equiv a_1\bmod q}$ and $c''_{D_2}=c'_{D_2}\mathbb 1_{D_2\equiv a_2\bmod q}$. By Lemma \ref{lemma_hyp}, the quantity above is bounded by $O\left( q^2\frac{T(\log T)}{ D' z^{1/4}}\right)$, so $S^\mathbf{A}(T;q,a)$ is bounded by \begin{equation*}
        \ll \frac{q^2 T\log T}{z^{1/4}}\sum_{\substack{D'\le T/ z^2}}\frac{\mu^2(D')\tau_{14}(D')}{D'4^{\omega(D')}}\ll\frac{q^2 T\log T}{z^{1/4}}\cdot (\log T)^{7/2}\,.\qedhere
    \end{equation*}
\end{proof}

\subsection{Type 3 error terms} We say that a dyadic multi-interval $\mathbf A$ is of \textbf{type 3} if it is not of type 2 or 1, if there exists an index $(i,j)$ such that $A_{ij}\ge C(T)$ and, for every $(k,l)$ linked to $(i,j)$, one has $A_{kl}\le z$ and at least one such linked $A_{kl}$ is not $1$. We further say that any  $\mcal D=\mcal L\mcal K$ belonging to such a multi-interval is a type 3 term if $K_{kl}>1$ for at least one index $(k,l)$ linked to $(i,j)$.
For a type 3 multi-interval $\mbf A$, define 
\begin{align*}G_h^{\mathbf A}(D)\coloneqq \mu^2(D)\sum_{\mcal D\in \D_{16}(D)}g(\mcal D)\mathbb 1_{\mathbf A, \text{type 3}}(\mcal D)\qquad \text{and}\qquad S_h^{\mathbf A}(T;q,a)\coloneqq \sum_{\substack{D\le T\\D\equiv a\bmod q\\D\sfree}}G_h^\mathbf A(D)\,.
\end{align*}

Setting $\hat K\coloneqq \prod_{(ij)\sim(kl)}K_{kl}$ and $\hat \ell\coloneqq \prod_{(ij)\sim(kl)}\ell_{kl}$; then $1<\hat \ell\hat K\ll z^{6}$. Also write $K=K'K_{ij}$ for some $K'$ coprime to $K_{ij}$ and, given $\mcal K\in \D_{16}(K)$, write $\mcal K= \mcal K'\times \{K_{ij}\}$ for some $\mcal K'\in \D_{15}(K')$. Then we can write $g$ as
$$g(\mcal D)=g(\mcal L)4^{-\omega( K)}\epsilon_{\mcal L,\mcal K'}\leg{D_{ij}}{\hat \ell\hat K}\psi_{\mcal L,\mcal K'}(D_{ij})\,,$$ 
where $\ord{\epsilon_{\mcal L,\mcal K'}}=1$, $\psi_{\mcal L,\mcal K'}$ is a character modulo 8, and both are independent of $D_{ij}$.
The character $\chi_{\mcal L,\mcal K'}\coloneqq \psi_{\mcal L,\mcal K'}\leg{\cdot}{\hat \ell\hat K}$ has modulus $8\hat \ell\hat K\ll z^{6}\ll (\log T)^{6B+N}$; hence we can employ a uniform version of Siegel--Walfisz to prove cancellation for type 3 multi-intervals.
\begin{lemma}\label{lemma6}
    Let $M>0$ and let $Q,r>0$. Let $\chi$ be a non-principal character modulo $Q$. Then there is a positive constant $c=c(M)$ such that
    $$\sum_{\substack{n\le x}}4^{-\omega(n)}\chi(n)\mu^2(nr)\ll_{M} x\tau(r)\exp\left(-c\sqrt{\log x}\right)\,,$$
    uniformly for $Q\le (\log x)^M$.
\end{lemma}

\begin{prop}\label{prop_error3}
    In the setting of Theorem \ref{teo_HB_arithm_prog_odd_restr}, fix $B,\eta>0$. Let $z(T), C(T)$ be as in (\ref{def_zC}), and let $\mathbf A$ be a type 3 multi-interval. Then there is a positive constant $c(N,B,\eta)$ such that $$\ord{S^\mathbf A_h(T;q,a)}\ll \frac{T(\log T)^{8}}{ e^{c(\log T)^{\eta/2}}}\,,$$ where the implied constant depends on $N,B,\eta$.
\end{prop}
\begin{proof}
    Let $(i,j)$ be the index for which $A_{ij}>C$, $1<A_{kl}\le z$ for every $(k,l)$ linked to $(i,j)$ and $\hat K\coloneqq \prod_{(i,j)\sim (k,l)}K_{kl}\neq 1$. Then 
    \begin{align*}
        \ord{S^{\mathbf{A}}_h(T;q,a)}
        &=
        \Bigg|\sum_{\mcal L\in \D_{16}(\ell)}{g(\mcal L)}\sum_{\substack{K\le T/\ell\\K\equiv a'\bmod q'\\\gcd(K,\ell)=1}}\frac{\mu^2(K)}{4^{\omega(K)}}\sum_{\substack{\mcal K\in \D_{16}(K)\\\mcal L\mcal K\in \mathbf{A}\,,\hat K\neq 1}}\epsilon_{\mcal L,\mcal K'}\,\chi_{\mcal L,\mcal K'}(K_{ij})\Bigg|
        \\&\le
         \sum_{\mcal L\in \D_{16}(\ell)}4^{-\omega(\ell)}\summ{K'\le T/\ell C'\\\mcal L'\mcal K'\in\mathbf A', \hat K\neq 1\\\gcd(K',\ell)=1}\frac{\mu^2(K')}{4^{\omega(K')}}\Bigg|\summ{C'\le n\le T/\ell K'\\n\in [A'_{ij},2A'_{ij})\\\gcd(n,\ell K')=1\\nK'\equiv a'\bmod q' }\mu^2(n)4^{-\omega(n)}\chi(n)\Bigg|\,,
    \end{align*}
    where we have used the triangle inequality together with $\ord{\epsilon_{\mcal L,\mcal K'}}\le 1$ and, for simplicity, we have relabelled $n=K_{ij}$, $C'=C/\ell_{ij}$ and $A_{ij}'=A_{ij}/\ell _{ij}$. 
    
    Letting $w\coloneqq \gcd(a',q')$, $q''\coloneqq q'/w$, and $a''\coloneqq a'/w$, observe that $\gcd(a,q)=w\ell$, and so $w$ must be square-free and coprime to $\ell$. Also observe that $\gcd(w,q'')=1$, otherwise $w$ would share a common factor with $\ell$. Set $\lambda\coloneqq \gcd(n,w)$, $n_1\coloneqq n/\lambda$, $\kappa\coloneqq \gcd(K',w)$, and $K_1'\coloneqq K'/\kappa$. Observe that the coprimality conditions in the sums force $\gcd(n_1K_1', \ell \lambda \kappa)=1$ and the square-free condition forces $\gcd(n_1,K_1')=1$. Now $$nK'\equiv a'\bmod q'\qquad \text{if and only if }\qquad\begin{cases}
        \lambda \kappa=w\\n_1K'_1\equiv a''\bmod q''\,.
    \end{cases}$$ 
    In particular, since $\gcd(a'',q'')=1$, we also need $\gcd(n_1K_1',q'')=1$. Part of the sum above can be rewritten as
    \begin{align*}
        \summ{w=\lambda\kappa}4^{-\omega(\lambda \kappa)}\summ{K_1'\le T/\ell C'\kappa\\\mcal L'\mcal K'\in\mathbf A',\,\hat K\neq 1\\\gcd(K'_1,\ell\lambda\kappa q'')=1}\frac{\mu^2(K_1')}{4^{\omega(K_1')}}\Bigg|\chi(\lambda)\summ{C'/\lambda\le n_1\le T/\ell K'\lambda\\n_1\lambda\in [A'_{ij},2A'_{ij})\\\gcd(n_1,\ell \lambda\kappa K'_1q'')=1\\n_1\equiv a''/K'_1\bmod q'' }\mu^2(n_1)4^{-\omega(n_1)}\chi(n_1)\Bigg|
    \end{align*}
    Using orthogonality of characters and the triangle inequality, we can bound the sum over $n_1$ in absolute value by 
    \begin{equation}\label{eq11}\ll \frac 1 {\phi(q'')}\summ{\theta\bmod q''}\Bigg|\overline{\theta}(a''/K_1')\summ{C'/\lambda\le n_1\le T/\ell K'\lambda\\n_1\lambda\in [A'_{ij},2A'_{ij})\\\gcd(n_1,\ell \lambda\kappa K'_1q'')=1 }\mu^2(n_1)4^{-\omega(n_1)}\chi(n_1)\theta(n_1)\Bigg|\,.\end{equation} Now $\theta\chi$ can never be the trivial character: indeed, since $\chi=\leg{\cdot}{\hat \ell \hat K}\psi$ with $\psi$ modulo $8$, if $\theta\chi$ were trivial then $\hat K$ would divide $q''$, but this cannot happen because $\hat K\mid K'=\kappa K_1'$, $\gcd(\kappa,q'')\mid \gcd(w,q'')=1$, $\gcd(K_1',q'')=1$, and $\hat K\neq 1$.
    Since $\theta\chi$ is a non-trivial character of modulus $O(q''\hat\ell\hat K)=O\left((\log T)^{N+6B}\right)$, Lemma \ref{lemma6} applies and (\ref{eq11}) can be bounded by
    $$\ll\frac {\tau(\ell\lambda \kappa K'_1q'')T}{\ell K'\lambda e^{c'\sqrt{\log {(A_{ij}/\ell_ij\lambda)}}}}\ll\frac {\tau(\ell\lambda \kappa K'_1q'')T}{\ell K'\lambda e^{c''(\log T)^{\eta/2}}}\,,$$ and therefore $\ord {S_h^\mathbf A(T;q,a)}$ is bounded by \begin{align*}
        &\ll \frac {\tau( q'')T}{ e^{c''(\log T)^{\eta/2}}}
        \sum_{\mcal L\in \D_{16}(\ell)}\frac{\tau(\ell)}{\ell4^{\omega(\ell)}}
        \summ{w=\lambda\kappa}\frac{\tau(\lambda)}{\lambda4^{\omega(\lambda)}}
        \summ{K'\le T/\ell C'\\\mcal L'\mcal K'\in\mathbf A',\, \hat K\neq 1\\\gcd(K'_1,\ell\lambda q'')=1\\\kappa\mid K'}\frac{\mu^2(K')\tau(K')}{K'4^{\omega(K')}}\\
        &\ll  \frac {\tau( q'')T}{ e^{c''(\log T)^{\eta/2}}}\frac{8^{\omega(\ell)}}{\ell}
           2^{\omega(w)} (\log T)^{15/2}\ll \frac{T(\log T)^8}{ e^{c''(\log T)^{\eta/2}}}\,.\qedhere
    \end{align*}    
    \end{proof}

\subsection{Type 4 error terms} Excluding type 1 terms means that there must be at least four variables larger than $C(T)$. Excluding type 2--3 terms implies that the large $D_{ij}$ are pairwise unlinked and that every $K_{kl}$ linked to some large $D_{ij}$ must be $1$. Lemma \ref{lemma_struct} shows that any set of at least four pairwise unlinked indices must consist of exactly four pairwise unlinked indices, and that every index outside this set must be linked to one of those four variables. We let $\mathscr B\subset \mcal I^4$ denote the set of 24 possible 4-tuples of pairwise unlinked indices, as identified in~\cite[Lemma 9]{HB}.

Given any set of indices $\mcal J\in \mathscr B$, we say that $\mcal M\in \N^{16}$ is supported on $\mcal J$ if $M_{kl}=1$ for every $(k,l)\notin \mcal J$. We say that $\mcal D=\mcal L\mcal K$ is a \textbf{type 4} term if $\mcal K$ is supported on $\mcal J$ and $D_{ij}\ge C(T)$ for every $(i,j)\in \mcal J$. Furthermore, we say that $D$ is of \textbf{type 4a} if $\mcal L$ is also supported on $\mcal J$, and of \textbf{type 4b} otherwise.
\begin{remark}
    Error terms belonging to a type 3 interval can be either of type 3 or 4b. Heath-Brown treated type 3 terms uniformly under~\cite[Lemma 6]{HB}, while restricting to arithmetic progressions imposes a different treatment of types 3 and 4b.
\end{remark}

To deal with type 4 terms, we highlight, when possible, two distinct characters in the formula for $g(\mcal L\mcal K)$ and obtain cancellation through the following version of Siegel--Walfisz~\cite[Lemma 10]{HB}.
\begin{lemma}\label{lemma10}
Let $M>0$, let $r,Q\in\Z_{>0}$, and let $\chi_1,\chi_2$ be two distinct characters modulo $Q$. Given $C_1,C_2\ge3$, call $C_0=\min\{C_1,C_2\}$. Then there is a positive constant $c=c(M)$ such that
$$
\sum_{\substack{mn\le x\\
m\ge C_1,\,n\ge C_2\\
}}
4^{-\omega(mn)}
\mu^2(mnr)\chi_1(m)\chi_2(n)
\ll_M
\tau(r)x\log x\exp\left(-c\sqrt{\log C_0}\right),
$$
uniformly for $x\ge C_1C_2$ and $Q\le(\log C_0)^M$.
\end{lemma}

For $\mcal J\in \mathscr B$, define \begin{align*}&G_h^{(\mcal J,\text{b})}(D)\coloneqq\mu^2(D) \summ{\mcal L\in \D_{16}(\ell)\\\mcal L\text{ not supp. on }\mcal J}\summ{\mcal K\in \D_{16}(K)\\\mcal K\text{ supp. on }\mcal J\\\mcal L\mcal K\text{ of type 4b}}g(\mcal L\mcal K)\qquad\text{and}\qquad  S_h^{(J,\text{b})}(T;q,a)\coloneqq \sum_{\substack{D\le T\\D\equiv a\bmod q\\D\sfree}}G_h^{(\mcal J,\text{b})} (D)\,.
\end{align*} By explicitly writing $g$ for type 4 terms supported on $\mcal J\in \mathscr B$, we have \begin{equation}\label{eq_gLK}
g(\mcal L\mcal K)=g(\mcal L)4^{-\omega(K)}\leg {-1}{\alpha(\mcal K)}\leg2{\beta(\mcal K)}\prod_{\substack{\mathbf{u}\in \mcal J}}\leg {K_{\mathbf u}}{\hat\ell_\mathbf u}(-1)^{\rho}\,,\end{equation}
where $$\hat\ell_{\mathbf u}=\prod_{(i,j)\text{ linked to }\mathbf u}\ell_{ij}\qquad\text{and}\qquad \rho=\summ{\{\mathbf u,\mathbf{v}\}\subset \mcal J\\\text{joined}}\frac{K_{\mathbf{u}}-1}{2}\frac{K_{\mathbf{v}}-1}{2}+\summ{\mathbf u\in \mcal J,(i,j)\in \mcal I\\\mathbf u\preceq (i,j)}\frac{K_{\mathbf{u}}-1}{2}\frac{\ell_{ij}-1}{2}\,.$$ In particular, $\rho$ is the term arising from quadratic reciprocity.

Now consider type 4b terms $\mcal D=\mcal L\mcal K$. By definition, there exists $(i,j)\notin \mcal J$ such that $\ell_{ij}>1$. By Lemma \ref{lemma_struct}, we can find $\mathbf u,\mathbf v\in \mcal J$ such that $(i,j)$ is linked to $\mathbf u$ but not to $\mathbf v$. This implies that $\ell_{ij}\mid \hat\ell_\mathbf u$ but $\ell_{ij}\nmid \hat\ell_{\mathbf v}$, and hence $\leg{\cdot}{\hat\ell_{\mathbf u}}$ and $\leg{\cdot}{\hat\ell_{\mathbf v}}$ are different characters. Let $K_{\mbf a},K_{\mbf b}$ denote the remaining two variables indexed by $\mcal J$. Rewriting (\ref{eq_gLK}) for type 4b with $\mcal L$ fixed, we obtain \begin{equation}\label{eq_gLK'}
    g(\mcal L\mcal K)=g(\mcal L)4^{-\omega(K)}\chi_\mbf u(K_\mbf u)\chi_\mbf v(K_\mbf v) \sigma(K_\mbf u,K_\mbf v)(-1)^{\rho'}\,,
\end{equation} where $\chi_{\mbf u},\chi_{\mbf v}$ are two different characters whose modulus is $\ll\ell\ll (\log T)^N/2$, $\rho'\in \Z$, and all three depend only on $\mcal L,K_{\mbf a},K_{\mbf b}$. Moreover,
\begin{equation}\label{def_sigma}\sigma(K_\mbf u,K_\mbf v)=(-1)^{\mathbb 1_{\mbf u,\mbf v\text{ joined}}\frac{K_\mbf u-1}{2}\frac{K_\mbf v-1}{2}}=\left[\frac 1 2\left(1+\leg {-1}{K_\mbf u}+\leg{-1}{K_\mbf v}-\leg {-1}{K_\mbf uK_\mbf v}\right)\right]^{\mathbb 1_{\mbf u,\mbf v\text{ joined}}}\,.\end{equation}

\begin{prop}\label{prop_err4b}
    In the setting of Theorem \ref{teo_HB_arithm_prog_odd_restr}, fix $\eta>0$ and $C(T)$ as in (\ref{def_zC}). Then, for every $\mcal J\in \mathscr B$, there is a positive constant $c(N,\eta)$ such that $$\ord{S^\mathbf {(\mcal J,\text{b})}_h(T;q,a)}\ll \frac{T(\log T)^2}{e^{c(\log T)^{\eta/2}}}\,,$$ where the implied constant depends on $N,\eta$.
\end{prop}
\begin{proof}
First, expand $\ord{S^{(\mcal J,\text{b})}_h}$ and bound it by 
$$\ord{S^\mathbf {(\mcal J,\text{b})}_h(T;q,a)}\le\summ{\mcal L\in \D_{16}(\ell)\\\mcal L\text{ not supp. on }\mcal J}4^{-\omega(\ell)}\Bigg|\summ{K\le T/\ell\\\gcd(K,\ell)=1\\K\equiv a'\bmod q'}\mu^2(K)\summ{\mcal K\in \D_{16}(K)\\\mcal L\mcal K \text{ of type 4}\\\mcal K \text{ supp. on }\mcal J}R(\mcal L,\mcal K)g(\mcal K)\Bigg|\,.$$ 
Fix a value of $\mcal L$ and, using the notation introduced above, bound the absolute value of the sum over $K$ as \begin{align*}
    \summ{K_\mbf aK_\mbf b\le T\ell_\mbf u\ell_\mbf v/\ell C^2\\C/\ell_\mbf a\le K_\mbf a,\,C/\ell_\mbf b\le K_\mbf b\\
    \gcd(K_\mbf a K_\mbf b,\ell)=1}{\mu^2(K_\mbf aK_\mbf b)}\Bigg|\summ{K_\mbf u K_\mbf v\le T/\ell K_\mbf a K_\mbf b\\C/\ell_\mbf u\le K_\mbf u,\,C/\ell_\mbf  vK_\mbf v\\\gcd(K_\mbf u K_\mbf v,K_\mbf aK_\mbf b\ell)=1\\K_\mbf uK_\mbf vK_\mbf aK_\mbf b\equiv a'\bmod q'}\mu^2(K_\mbf uK_\mbf v)R(\mcal L\mcal K)g(\mcal K)\Bigg|
\end{align*} Set $K'=K_\mbf aK_\mbf b$. Then, by (\ref{eq_gLK'}), we can bound this by $$ \le \summ{K_\mbf aK_\mbf b\le T\ell_\mbf u\ell_\mbf v/\ell C^2\\C/\ell_\mbf a\le K_\mbf a,\,C/\ell_\mbf b K_\mbf b\\
    \gcd(K_\mbf a K_\mbf b,\ell)=1}{\mu^2(K_\mbf aK_\mbf b)}\Bigg|\summ{K_\mbf u K_\mbf v\le T/\ell K'\\C/\ell_\mbf  u\le K_\mbf u,\, C/\ell_\mbf  v\le K_\mbf v\\\gcd(K_\mbf u K_\mbf v,K'\ell)=1\\K_\mbf uK_\mbf vK'\equiv a'\bmod q'}\mu^2(K_\mbf uK_\mbf v)\chi_\mbf u(K_\mbf u)\chi_\mbf v(K_\mbf v)\sigma(K_\mbf u,K_\mbf v)\Bigg|\,.$$
For simplicity, we assume that $\mbf u$ and $\mbf v$ are not joined, so that $\sigma(K_\mbf u,K_\mbf v)=1$; if $\mbf u$ and $\mbf v$ are joined, the same process applies after splitting the sum above into four character sums using (\ref{def_sigma}). Set $w=\gcd(a',q')$, $q''=q'/w$, $a''=a'/w$, $\kappa_\mbf u=\gcd(w,K_\mbf u)$, $K_\mbf u'=K_{\mbf u}/\kappa_\mbf u$, and use the same notation for $\mbf v,\mbf a,\mbf b$. Then the sum above is bounded by
$$ \summ{w=\kappa_\mbf a\kappa_\mbf b\kappa_\mbf u\kappa_\mbf v}\!\!\!\!\!\! 4^{-\omega(\kappa_\mbf u\kappa_\mbf v)}\!\!\!\!\!\!\!\!\!
\summ{K_\mbf aK_\mbf b\le T\ell_\mbf u\ell_\mbf v/\ell C^2\\C/\ell_\mbf  a\le K_\mbf a,\,C/\ell_\mbf  b\le K_\mbf b\\\kappa_\mbf a\mid K_\mbf a,\, \kappa_\mbf b\mid K_\mbf b\\    \gcd(K_\mbf a'K_\mbf b',\ell w)q''=1}\!\!\frac{\mu^2(K_\mbf aK_\mbf b)}{4^{\omega(K_\mbf a K_\mbf b)}}\Bigg|\!\!\!
\summ{K_\mbf u 'K_\mbf v'\le T/\ell K'\kappa_\mbf u\kappa_\mbf v\\C/\kappa_\mbf u\ell_\mbf u\le K_\mbf u',\,C/\kappa_\mbf v\ell_\mbf v\le K_\mbf v'\\\gcd(K_\mbf u' K_\mbf v',K'\ell w)=1\\K_\mbf u'K_\mbf v'\equiv a''(K'_\mbf aK_\mbf b')^{-1}\bmod q''}\!\!\!\!\!\!\!\!\mu^2(K_\mbf u'K_\mbf v')4^{-\omega(K_\mbf u' K_\mbf v')}\chi_\mbf u(K_\mbf u')\chi_\mbf v(K_\mbf v')\Bigg|\,.$$
By character orthogonality, we can bound the sum inside the absolute value by
$$\frac{1}{\phi(q'')}\summ{\theta\bmod q''}\Bigg|
\summ{K_\mbf u 'K_\mbf v'\le T/\ell K'\kappa_\mbf u\kappa_\mbf v\\C/\kappa_\mbf u\ell_\mbf u\le K_\mbf u',\,C/\kappa_\mbf v\ell_\mbf v\le K_\mbf v'\\\gcd(K_\mbf u' K_\mbf v',K'\ell w)=1}\mu^2(K_\mbf u'K_\mbf v')4^{-\omega(K_\mbf u' K_\mbf v')}[\theta\chi_\mbf u](K_\mbf u')[\theta\chi_\mbf v](K_\mbf v')\Bigg|\,.$$
As $\theta \chi_\mbf u$ and $\theta\chi_\mbf v$ are different characters whose modulus is at most $8\ell q''\ll (\log T)^N$, we can apply Lemma \ref{lemma10}, and, for fixed $\theta$, each sum is bounded by $\ll \frac{\tau(\ell K'\kappa_\mbf u\kappa_\mbf v)T(
\log T)}{\ell K'\kappa_\mbf u\kappa_\mbf v e^{c'(\log T)^{\eta/2}}}$. 
Combining this with the preceding estimates, $\ord {S_h^{(\mcal J,\text{b})}(T;q,a)}$ can be bounded by 
\begin{align*}
     \ll& \frac{T(
\log T)}{ e^{c'(\log T)^{\eta/2}}}\summ{\mcal L\in \D_{16}(\ell)\\\mcal L\text{ not supp. on }\mcal J}\frac{\tau(\ell)}{\ell4^{\omega(\ell)}}\summ{w=\kappa_\mbf a\kappa_\mbf b\kappa_\mbf u\kappa_\mbf v} \frac{\tau(w)}{w4^{\omega(w)}} 
\summ{K_\mbf a'K_\mbf b'\le T\ell_\mbf u\ell_\mbf v/\ell C^2\kappa_\mbf a\kappa_\mbf b\\C/\kappa_\mbf a \ell\mbf a\le K_\mbf a',\,C/\kappa_\mbf b\ell_\mbf b\le K_\mbf b'\\    \gcd(K_\mbf a'K_\mbf b',\ell w)=1}\frac{\mu^2(K_\mbf a'K_\mbf b')\tau(K_\mbf a' K_\mbf b')}{K_\mbf a 'K_\mbf b'4^{\omega(K_\mbf a' K_\mbf b')}}
\\\ll & \frac{T(
\log T)}{ e^{c'(\log T)^{\eta/2}}}\frac{8^{\omega(\ell)}}{\ell}\frac{2^{\omega(w)}}{w}(\log T)\ll \frac{T(
\log T)^2}{ e^{c'(\log T)^{\eta/2}}}\,. \qedhere
\end{align*}
\end{proof}
\begin{remark}
    Until now, we have been exploiting only the combinatorics of linked variables to search for some source of cancellation in the expression of $g$. From now on, one needs to consider also the contribution of the two characters $\leg {-1}\alpha$ and $\leg 2 \beta$. 
\end{remark}

We are now left with type 4a error terms. Since both $\mcal L$ and $\mcal K$ are supported on a fixed $\mcal J\in \mathscr B$, we can reunite $\mcal L\mcal K$ into a single $\mcal D$ supported on $\mcal J$.

Depending on the shape that $g$ takes on type 4a variables supported on $\mcal J\in \mathscr B$, we split $\mathscr B$ into
$$\mathscr B_1\coloneqq\{(10,20,30,40),(40,41,42,43)\},\qquad\qquad \mathscr B_{1/2}\coloneqq\{(20,12,32,42),(30,13,23,43)\}\,,$$
and $\mathscr B_0\coloneqq\mathscr B\setminus \mathscr B_1\setminus \mathscr B_{1/2}\,.$
In particular, Heath-Brown~\cite[\S\S 4--5]{HB} showed that if $\mcal D$ is a type 4a variable supported on $\mcal J$, then:\begin{itemize}
    \item for every $\mcal J\in \mathscr B_0$, one can relabel the four non-one variables as $D_{\mbf a},D_\mbf b,D_\mbf u,D_\mbf v$ so that
    $$g(\mcal D)=4^{-\omega(D)}(-1)^{\rho}\chi_\mbf u(D_\mbf u)\chi_\mbf v(D_\mbf v)\,,$$
    where $\rho\in \Z$, $\chi_\mbf u,\chi_\mbf v$ are two distinct characters modulo $8$, and $\rho, \chi_\mbf u,\chi_\mbf v$ depend only on $D_\mbf a,D_\mbf b$;
\item for every  $\mcal J\in \mathscr B_1$ $$g(\mcal D)=4^{-\omega(D)}\,;$$
\item for every $\mcal J\in \mathscr B_{1/2}$, after relabelling the variables as $D_{\mbf a},D_\mbf b,D_\mbf u,D_\mbf v$, $g$ rewrites as $$g(\mcal D)=\frac{4^{-\omega(D)}} 2\left[1+\leg{-1}{D_\mbf b D_\mbf u}+\leg{-1}{D_\mbf bD_\mbf v}-\leg{-1}{D_\mbf uD_\mbf v}\right]\,.$$
\end{itemize}    
The extra terms in $g$ involving two different characters can be cancelled using Lemma \ref{lemma10}, while the four choices of $\mcal J$ in $\mathscr B_1\cup \mathscr B_{1/2}$ produce, respectively, $1,1,1/2,1/2$ as a main term.

Define $$G^{(\mcal J,\text{a})}(D)=\summ{\mcal D\in  \D_{16}(D)\\\mcal D\text{ of type 4a}\\\mcal D\text{ supp. on }\mcal J}\left(g(\mcal D)-4^{-\omega(D)}\varsigma(\mcal J)\right)\qquad \text{with }\qquad \varsigma(\mcal I')=\begin{cases}
    1&\mcal I'\in \mathscr B_1\\
    1/2&\mcal I'\in \mathscr B_{1/2}\\
    0&\mcal I'\in \mathscr B_0
\end{cases}\,,$$ and define $S^{(\mcal J,\text{a})}(T,q;a)$ in the usual way.
\begin{prop}\label{prop_err4a}
    In the setting of Theorem \ref{teo_HB_arithm_prog_odd_restr}, fix $\eta>0$ and $C(T)$ be as in (\ref{def_zC}). Then, for every $\mcal J\in \mathscr B$, there is a positive constant $c(N,\eta)$ such that $$\ord{S^\mathbf {(\mcal J,\text{a})}_h(T;q,a)}\ll  \frac{T(\log T)^2}{e^{c(\log T)^{\eta/2}}}\,,$$ where the implied constant depends on $N,\eta$.
\end{prop}
\begin{proof}
    The argument follows the same lines as the proof of Proposition \ref{prop_err4b}: depending on the specific $\mcal J$, identify two variables $D_\mbf u,D_\mbf v$  whose characters appearing in the expression of $g(\mcal D)-4^{-\omega(D)}\varsigma(\mcal J)$ are different. After factoring out $\gcd(a,q)$, deal with the congruence constraint by using characters' orthogonality; then the two emerging characters remain different, and the upper bound follows by Lemma \ref{lemma10}. 
\end{proof}

Finally, after removing type 1--4 error terms, we are left with
 \begin{align*}
     \summ{D\le T\\D\equiv a \bmod q\\D\sfree}\sum_{\mcal J\in \mathscr B}\summ{\mcal D\in \D_{16}(D)\\\mcal D\text{ of type 4a}\\\mcal D\text{ supp. on }\mcal J}4^{-\omega(D)}\varsigma(\mcal J)&= 3 \summ{D\le T\\D\equiv a \bmod q\\D\sfree}4^{-\omega(D)}\#\{\mcal D'\in \D_4(D)| \forall \mbf i,\, D'_\mbf i\ge C \}\\  
     \end{align*}
The set in the sum above has cardinality $4^{-\omega(D)}-\#\{\mcal D'\in \mathbb D_{4}(D)\mid \exists \mbf i\, D'_\mbf i\le C\}$, so the quantity above produces the main term of Theorem \ref{teo_HB_arithm_prog_odd_restr} plus one final residual error term. Define
     $$G^{(\text{R})}(D)\coloneqq 3\mu^2(D)\#\{\mcal D'\in \D_4(D)| \exists \mbf i,\, D'_\mbf i< C \}$$
    and define $S^{(\text{R})}(T;q,a)$ in the usual way. Then
    \begin{equation}\label{eq_12}G^{(\text{R})}(D)\ll \mu^2(D)\summ{D'\le C\\D'\mid D}\leg 3 4 ^{\omega\left(D/D'\right)}\,.\end{equation}
    \begin{prop}\label{prop_errR}
 In the setting of Theorem \ref{teo_HB_arithm_prog_odd_restr}, fix $0<\eta<1/4$ and $C(T)$ as in (\ref{def_zC}). Then  $$\ord{S^\mathbf {(\text{R})}(T;q,a)}\ll \frac{T}{q(\log T)^{1/4-2\eta}}\,,$$ where the implied constant depends on $N,\eta$.
 \end{prop}
 \begin{proof}
    The proof exploits the same underlying mechanism as in the proof of Proposition \ref{prop_error1}, where the role of inequality (\ref{bound_g3}) is replaced by inequality (\ref{eq_12}).
 \end{proof}
    
    \begin{proof}[Proof of Theorem \ref{teo_HB_arithm_prog_odd_restr}.]
    Recapping our error term analysis, once $h$ is fixed, we have
    \begin{equation}2^{s(D)}=3+\Bigg[G^{(1)}+\summ{\mbf A \text{ multi-interval}\\\text{of type 2}}G^{\mbf A}+\summ{\mbf A \text{ multi-interval}\\\text{of type  3}}G_h^{\mbf A}+\sum_{\mcal J\in \mathscr B}\left(G^{(\mcal J,\text{a})}+G_h^{(\mcal J,\text{b})}\right)-G^{(\text{R})}\Bigg](D)\,.\end{equation}
    Recall that in the two sums above there are $\ll (\log T)^{16}$ relevant dyadic multi-intervals.
    Having chosen a residue class $a$ modulo $q$, we use the triangle inequality to separate the contribution of the different terms in $S(T;q,a)$, and we bound each term individually using Propositions \ref{prop_error2}, \ref{prop_error3}, \ref{prop_error1}, \ref{prop_err4a}, \ref{prop_err4b}, and \ref{prop_errR}. After choosing $B>0$ large enough in $z(T)$ and $\eta>0$ small enough in $C(T)$ (see (\ref{def_zC})), the claimed error term follows.
    \end{proof}
\subsection{2-Selmer on even arithmetic progressions.}
Here we prove a version of Theorem \ref{teo_HB_arithm_prog_odd_restr} for even arithmetic progressions. As we showed in \S\ref{sec_prel}, the 2-descent counting formula (\ref{2desc_even}) is structurally different for even twists; however, the functions $g$ in (\ref{defg}) and $g_2$ in (\ref{defg_even}) differ only in the factors $(\alpha,\beta)$ and $(\alpha_2,\beta_2)$. Hence we can repeat part of the proof of Theorem \ref{teo_HB_arithm_prog_odd_restr} verbatim for even arithmetic progressions. 

\begin{theorem}
    \label{teo_HB_arithm_prog_even_restr} Fix $N,\epsilon>0$. Let $T>3$ and $h\in\Z_{>0}$ with $h\le(\log T)^{N/2}$. Then, for every $q,a\in\Z_{>0}$ with $q,a\le (\log T)^N$, $q$ even, $a$ odd, $h^2\mid q$ and $\mu^2(q/h^2)=1$, we have $$S(T;q,2a)=3\#\left\{2D\le T\bigg|\,\begin{aligned}
        &2D\equiv 2a\bmod q\\[-1.5mm]&2D\sfree
    \end{aligned}\right\}+O\left(\frac T{q\,(\log T)^{1/4-\epsilon}}\right)\,,$$
    where the implied constant depends only on $N,\epsilon$.
\end{theorem}
\begin{proof}
In the proof of Theorem \ref{teo_HB_arithm_prog_odd_restr}, all the error-term analysis up to and including type 4b depends only on the links between the indices in $\mcal I$. Hence, one can analogously define $G^{\mbf A .2},G_h^{\mbf A.2}, G^{(1).2}, G_{h}^{(\mcal J,\text{b}).2}$ by using $g_2$ in place of $g$ in the definitions of $G^{\mbf A},G_h^{\mbf A}, G^{(1)}, G_{h}^{(\mcal J,\text{b})}$, and obtain the corresponding versions of Propositions \ref{prop_error2}, \ref{prop_error3}, \ref{prop_error1}, and \ref{prop_err4b} for even twists.
The influence of $\alpha_2$ and $\beta_2$ emerges only later in the proof, when we identify which $\mcal J\in \mathscr B$ contribute to the main term. To deal with type 4a terms supported on $\mcal J$, we identify two indices $\mbf u,\mbf v\in \mcal J$ whose respective characters $\chi_\mbf u,\chi_\mbf v$ in the expression of $g_2$ are different; this yields cancellation via Lemma \ref{lemma10}.

Observe that $\leg 2\cdot$ is a primitive character modulo 8, and recall that in (\ref{defg_even}) we have $\beta_2(\mcal D)    =D_{13}D_{23}D_{31}D_{32}D_{34}D_{43}$. Whenever there are $\mbf u,\mbf v\in \mcal J$ such that $\mbf u$ belongs to the list of variables appearing in $\beta_2$, while $\mbf v$ does not, then $\chi_\mbf u$ is a primitive character modulo 8 and $\chi_\mbf v$ is a character modulo 4; in particular, $\chi_\mbf u$ and $\chi_\mbf v$ cannot coincide. Using Table (\ref{table_indexes}), one sees that this happens for every $\mcal J\in \mathscr B_0'\setminus \{(10,20,30,40)\}$, where we set
$$\mathscr B'_1\coloneqq \{(10,20,21,12),(10,41,14,40),(24,20,42,40)\}\qquad\text{and}\qquad \mathscr B'_0\coloneqq \mathscr B\setminus \mathscr B'_1\,.$$
Then, if $\mcal D$ is a type 4a variable supported on $\mcal J$, we have: \begin{itemize}
    \item for every  $\mcal J\in \mathscr B_1'$, $$g_2(\mcal D)=4^{-\omega(D)}\,;$$
    \item for   $\mcal J=(10,20,30,40)$, $$g_2(\mcal D)=4^{-\omega(D)}\leg {-1}{D_{30}}\,;$$
    \item for every  $\mcal J\in \mathscr B_0'\setminus \{(10,20,30,40)\}$, after relabelling the variables, $$g_2(\mcal D)=4^{-\omega(D)}(-1)^{\rho}\chi_\mbf u(D_\mbf u)\chi_\mbf v(D_\mbf v)\sigma(D_\mbf u,D_\mbf v)\,,$$ where $\chi_\mbf u,\chi_\mbf v, \rho$ are independent of $D_{\mbf u},D_\mbf v$, $\chi_\mbf u$ is a primitive character modulo 8, $\chi_\mbf v$ is a character modulo 4, and $\sigma $ is a sum of characters modulo 4 as defined in (\ref{def_sigma}).
\end{itemize}    
Without giving further details, the proof continues as in Theorem \ref{teo_HB_arithm_prog_odd_restr}, where the identification of different characters allows cancellation for type 4a terms supported on $\mcal J\in \mathscr B'_0$, while the main term $3$ is produced by those type 4a terms supported on the three four-tuples in $\mathscr B_1'$.
\end{proof}

\subsection{A counterexample over short arithmetic progressions}
A recent breakthrough by OpenAI on the quasi-Riemann hypothesis~\cite{OAI} claims a zero-free region for every Dirichlet $L$-function. As a consequence, a stronger version of the Siegel--Walfisz method is available, and this would strengthen our application of Lemmas \ref{lemma6} and 
\ref{lemma10}. It is then natural to try to extend Theorem \ref{teo_HB_arithm_prog} to every $q\le T^\gamma$ for some small $\gamma>0$. As the following example shows, this is not possible without extra assumptions on $q$ and $a$. 

\begin{example}\label{ex_sap}
    Let $m>0$. Heath-Brown~\cite[\S 7]{HB} constructs $m$ primes $p_1,\ldots, p_m$ such that $p_i\equiv 1\bmod 8$ for every $1\le i\le m$,$\leg {p_i}{p_j}=1$ for every $1\le i\neq j\le m$, and $m\gg \sqrt{\log h}$, where $h\coloneqq p_1\cdot\ldots\cdot p_m$. Set  $a=h,$ $q=8h^2,$ and $T=h^A$ for some $A>2$. Then $$S(T;8h^2,h)=\summ{D\le T\\D\equiv h\bmod 8h^2\\D\sfree}2^{s(D)}\ge \summ{p<T/h\\p\equiv 1\bmod 8h\\p\text{ prime}}2^{s(ph)}\,.$$ The restriction $p\equiv 1\bmod 8h$ implies $\leg {p}{p_i}=1$ for every $1\le i\le m$, and so formulae (\ref{defg},\ref{2desc}) imply that $$g(ph)=4^{-m-1}\qquad\text{and}\qquad2^{s(ph)}=4^{m+1}\,.$$
By QRH, for $A$ large enough we have $$\#\{p\le T/h \text{ prime}: p\equiv 1\bmod 8h\}\asymp \frac{T}{h\phi(8h)\log T}\,,$$ and so 
$$S(T;q,a)\gg \frac{T4^m}{\phi(8h^2)\log T}\gg \frac T q\frac{4^m} {\log T}.$$ Now $m\gg\sqrt{\log h}$ implies $4^m/\log T\to \infty$ as $m\to \infty$, and so Theorem \ref{teo_HB_arithm_prog} cannot hold for such a choice of $q$ and $a$.
\end{example}

\begin{remark}
    In propositions \ref{prop_error3} and  \ref{prop_err4a} the contribution of $\mcal L$ can be easily absorbed as $\ell$ is small compared to $T$. Example \ref{ex_sap} shows that this is no longer possible when $\ell$ is sufficiently large; for instance, it can fail for $\ell\gg T^{\frac{(\log \log T)^3}{\log T}}$.
    \end{remark}



\bibliographystyle{plainurl}
\bibliography{sample}

@article{HB,
  AUTHOR   = {Heath-Brown, D. R.},
  TITLE    = {The size of {Selmer} groups for the congruent number problem},
  JOURNAL  = {Invent. Math.},
  FJOURNAL = {Inventiones Mathematicae},
  VOLUME   = {111},
  YEAR     = {1993},
  NUMBER   = {1},
  PAGES    = {171--195},
  DOI      = {10.1007/BF01231285},
  URL      = {https://doi.org/10.1007/BF01231285},
}

@article{Wi,
  AUTHOR   = {Wilson, C.},
  TITLE    = {General bilinear forms in the {Jacobi} symbol over hyperbolic regions},
  JOURNAL  = {Monatsh. Math.},
  FJOURNAL = {Monatshefte f{\"u}r Mathematik},
  VOLUME   = {202},
  YEAR     = {2023},
  NUMBER   = {2},
  PAGES    = {435--451},
  DOI      = {10.1007/s00605-023-01848-9},
  URL      = {https://doi.org/10.1007/s00605-023-01848-9},
}

@article{Na,
  AUTHOR   = {Nair, M.},
  TITLE    = {Multiplicative functions of polynomial values in short intervals},
  JOURNAL  = {Acta Arith.},
  FJOURNAL = {Acta Arithmetica},
  VOLUME   = {62},
  YEAR     = {1992},
  NUMBER   = {3},
  PAGES    = {257--269},
  DOI      = {10.4064/aa-62-3-257-269},
  URL      = {https://doi.org/10.4064/aa-62-3-257-269},
}

@article{HB2,
  AUTHOR   = {Heath-Brown, D. R.},
  TITLE    = {The size of {Selmer} groups for the congruent number problem. {II}},
  JOURNAL  = {Invent. Math.},
  FJOURNAL = {Inventiones Mathematicae},
  VOLUME   = {118},
  YEAR     = {1994},
  NUMBER   = {2},
  PAGES    = {331--370},
  NOTE     = {With an appendix by P. Monsky},
  DOI      = {10.1007/BF01231536},
  URL      = {https://doi.org/10.1007/BF01231536},
}

@article{SD,
  AUTHOR   = {Swinnerton-Dyer, P.},
  TITLE    = {The effect of twisting on the $2$-{Selmer} group},
  JOURNAL  = {Math. Proc. Cambridge Philos. Soc.},
  FJOURNAL = {Mathematical Proceedings of the Cambridge Philosophical Society},
  VOLUME   = {145},
  YEAR     = {2008},
  NUMBER   = {3},
  PAGES    = {513--526},
  DOI      = {10.1017/S0305004108001588},
  URL      = {https://doi.org/10.1017/S0305004108001588},
}

@article{Ka,
  AUTHOR   = {Kane, D. M.},
  TITLE    = {On the ranks of the $2$-{Selmer} groups of twists of a given
              elliptic curve},
  JOURNAL  = {Algebra Number Theory},
  FJOURNAL = {Algebra \& Number Theory},
  VOLUME   = {7},
  YEAR     = {2013},
  NUMBER   = {5},
  PAGES    = {1253--1279},
  DOI      = {10.2140/ant.2013.7.1253},
  URL      = {https://doi.org/10.2140/ant.2013.7.1253},
}

@article{Ro,
  AUTHOR   = {Rhoades, R. C.},
  TITLE    = {$2$-{Selmer} groups and the {Birch--Swinnerton-Dyer} conjecture
              for the congruent number curves},
  JOURNAL  = {J. Number Theory},
  FJOURNAL = {Journal of Number Theory},
  VOLUME   = {129},
  YEAR     = {2009},
  NUMBER   = {6},
  PAGES    = {1379--1391},
  DOI      = {10.1016/j.jnt.2009.01.015},
  URL      = {https://doi.org/10.1016/j.jnt.2009.01.015},
}

@misc{Sm,
  AUTHOR        = {Smith, A.},
  TITLE         = {The {Birch} and {Swinnerton-Dyer} conjecture implies
                   {Goldfeld}'s conjecture},
  YEAR          = {2025},
  NOTE          = {Preprint},
  EPRINT        = {2503.17619},
  ARCHIVEPREFIX = {arXiv},
  PRIMARYCLASS  = {math.NT},
  URL           = {https://arxiv.org/abs/2503.17619},
}

@article{BS,
  AUTHOR   = {Bhargava, M. and Shankar, A.},
  TITLE    = {Binary quartic forms having bounded invariants, and the
              boundedness of the average rank of elliptic curves},
  JOURNAL  = {Ann. of Math. (2)},
  FJOURNAL = {Annals of Mathematics. Second Series},
  VOLUME   = {181},
  YEAR     = {2015},
  NUMBER   = {1},
  PAGES    = {191--242},
  DOI      = {10.4007/annals.2015.181.1.3},
  URL      = {https://doi.org/10.4007/annals.2015.181.1.3},
}

@book{Si,
  AUTHOR    = {Silverman, J. H.},
  TITLE     = {The arithmetic of elliptic curves},
  SERIES    = {Graduate Texts in Mathematics},
  VOLUME    = {106},
  EDITION   = {Second},
  PUBLISHER = {Springer},
  ADDRESS   = {Dordrecht},
  YEAR      = {2009},
  PAGES     = {xx+513},
  DOI       = {10.1007/978-0-387-09494-6},
  URL       = {https://doi.org/10.1007/978-0-387-09494-6},
}

@article{KMR,
  AUTHOR   = {Klagsbrun, Z. and Mazur, B. and Rubin, K.},
  TITLE    = {A {Markov} model for {Selmer} ranks in families of twists},
  JOURNAL  = {Compos. Math.},
  FJOURNAL = {Compositio Mathematica},
  VOLUME   = {150},
  YEAR     = {2014},
  NUMBER   = {7},
  PAGES    = {1077--1106},
  DOI      = {10.1112/S0010437X13007896},
  URL      = {https://doi.org/10.1112/S0010437X13007896},
}

@article{Ch,
  AUTHOR   = {Chan, S.},
  TITLE    = {The average number of integral points on the congruent number curves},
  JOURNAL  = {Adv. Math.},
  FJOURNAL = {Advances in Mathematics},
  VOLUME   = {457},
  YEAR     = {2024},
  PAGES    = {Paper No. 109946, 31},
  DOI      = {10.1016/j.aim.2024.109946},
  URL      = {https://doi.org/10.1016/j.aim.2024.109946},
}

@misc{OAI,
  author = {{OpenAI}},
  title = {{The quasi-{R}iemann hypothesis:
           a zero-free Half-Plane
           $\mathrm{Re}(s)>7/8$}},
  year = {2026},
  note = {OpenAI Math Release preprint,
          30 September 2026},
  url = {https://github.com/openai/math/blob/main/preprints/The-Quasi-Riemann-Hypothesis-September-30-2026/paper.pdf}
}

@article{FK,
  AUTHOR   = {Fouvry, {\'E}. and Kl{\"u}ners, J.},
  TITLE    = {On the $4$-rank of class groups of quadratic number fields},
  JOURNAL  = {Invent. Math.},
  FJOURNAL = {Inventiones Mathematicae},
  VOLUME   = {167},
  YEAR     = {2007},
  NUMBER   = {3},
  PAGES    = {455--513},
  DOI      = {10.1007/s00222-006-0021-2},
  URL      = {https://doi.org/10.1007/s00222-006-0021-2},
}

@mastersthesis{Mi,
  AUTHOR  = {Milovic, D.},
  TITLE   = {On the $4$-rank of class groups of quadratic number fields},
  SCHOOL  = {Universit{\'e} Paris-Sud 11},
  ADDRESS = {Orsay},
  YEAR    = {2013},
  TYPE    = {Master's thesis},
  URL     = {https://algant.eu/documents/theses/milovic.pdf},
}

\end{document}